\documentclass[11 pt, twoside]{amsart}
\usepackage{parskip}
\usepackage{amsthm,amsmath,amssymb,oldgerm}
\usepackage[a4paper, margin=2.5cm]{geometry}
\usepackage{mathrsfs}
\usepackage{verbatim}
\usepackage[all]{xy}
\theoremstyle{plain}
\newtheorem{thm}{Theorem}[section]
\newtheorem*{theorem*}{Main theorem}

\theoremstyle{definition}
 
\theoremstyle{definition}
\newtheorem{rem}[thm]{Remark}

\let\H\relax

\let\O\relax
\let\dh\relax

\DeclareMathOperator{\bkx} {\mathcal{B}_{\Omega_{\it{X}}}^{\it{k}}}

\DeclareMathOperator{\H}{\mathbb{H}}
  
\DeclareMathOperator{\G}{\Gamma}

\DeclareMathOperator{\R}{\mathbb{R}} 
 
\DeclareMathOperator{\O}{\Omega_{{\it{X}}}} 
\DeclareMathOperator{\Ok}{\Omega_{{\it{X}}}^{\otimes{\it{k}}}}

\DeclareMathOperator{\rx}{{\it{r}}_{{\it{X}}}} 
\DeclareMathOperator{\dh}{{\it{d}}_{\mathrm{hyp}}}
\title[Estimates of the Bergman kernel]{Estimates of the Bergman kernel on a hyperbolic Riemann surface of finite volume-II}
\author{Anilatmaja Aryasomayajula}
\address{Department of Mathematics, Indian Institute of Science Education and Research Tirupati, 
Tirupati-517507, India.}
\email{anil.arya@iisertirupati.ac.in}
\author{Priyanka Majumder}
\address{Department of Mathematics, Indian Institute of Science Education and Research Tirupati, 
 Tirupati-517507, India.}
\email{pmpriyanka57@gmail.com }
\begin{document}
%%%%%%%%%%%%%%%%%%%%%%%%%%%%%%%%%%%%%%%%%%%%%%%%%%%%%%%%%%%%%%%%%%%%%%%%%%%%%%%%%%%%%%%%%%%%%
\begin{abstract}
\vspace{0.2cm}\noindent
In this article, we derive off-diagonal estimates of the Bergman kernel associated to the tensor-powers of the cotangent bundle defined on a hyperbolic Riemann surface of finite volume, when the distance between the points is less than injectivity radius. We then use these estimates to derive estimates of the Bergman kernel along the diagonal.

\vspace{0.2cm}\noindent
R{\tiny{\'ESUME}}. Dans cet article, nous dérivons des estimations non-diagonales du noyau de Bergman associé aux puissances tensorielles du faisceau cotangent défini sur une surface de Riemann hyperbolique de volume fini, lorsque la distance entre les points est inférieure au rayon d'injectivité. Nous utilisons ensuite ces estimations pour dériver des estimations du noyau de Bergman le long de la diagonale.
%%%%%%%%%%%%%%%%%%%%%%%%%%%%%%%%%%%%%%%%%%%%%%%%%%%%%%%%%%%%%%%%%%%%%%%%%%%%%%%%%%%%%%%%%%%%%

\vspace{0.15cm}\noindent
{\em Mathematics Subject Classification (2010)}: 32A25, 30F30, 30F35. {\em Keywords}: Bergman kernels. 
\end{abstract}
%%%%%%%%%%%%%%%%%%%%%%%%%%%%%%%%%%%%%%%%%%%%%%%%%%%%%%%%%%%%%%%%%%%%%%%%%%%%%%%%%%%%%%%%%%%%%
\maketitle
\vspace{-1cm}
%%%%%%%%%%%%%%%%%%%%%%%%%%%%%%%%%%%%%%%%%%%%%%%%%%%%%%%%%%%%%%%%%%%%%%%%%%%%%%%%%%%%%%%%%%%%%
\section{Introduction}
In \cite{am}, using elementary methods from complex analysis, we derive off-diagonal estimates of the Bergman kernel associated to tensor-powers of the cotangent bundle defined on a hyperbolic Riemann surface of finite volume, when the distance between the points is greater than injectivity radius, in both the compact and noncompact setting. In this article, refining the arguments from \cite{am}, we derive off-diagonal estimates of the Bergman kernel, when the distance between the points is less than  injectivity radius. This article complements the results from \cite{am}. 
%%%%%%%%%%%%%%%%%%%%%%%%%%%%%%%%%%%%%%%%%%%%%%%%%%%%%%%%%%%%%%%%%%%%%%%%%%%%%%%%%%%%%%%%%%%%%

\vspace{0.1cm}
Several mathematicians including the likes of Tian, Zelditch, Ma, Marinsecu et al. have derived estimates of Bergman kernels associated to high tensor-powers of  line bundles defined over complex manifolds. We refer the reader to the introduction in \cite{am}, for an elaborate discussion on similar results from literature. However, we mention the results from \cite{au-ma}, \cite{au-ma2}, and \cite{am}, which have led to the culmination of this article. 
%%%%%%%%%%%%%%%%%%%%%%%%%%%%%%%%%%%%%%%%%%%%%%%%%%%%%%%%%%%%%%%%%%%%%%%%%%%%%%%%%%%%%%%%%%%%%

\vspace{0.1cm}
Let $X$ be a noncompact Riemann surface, whose natural metric has singularities of Poincar\'e type at a finite set. Let $\mathcal{L}$ be a holomorphic line bundle whose curvature form is a scalar multiple of the hyperbolic metric outside a compact subset of $X$. In  \cite{au-ma}, Auvray, Ma, and Marinescu have derived optimal estimates of  $\mathcal{C}^{n}$-norms of the Bergman kernel associated to tensor-powers of $\mathcal{L}$, along the diagonal. 
%%%%%%%%%%%%%%%%%%%%%%%%%%%%%%%%%%%%%%%%%%%%%%%%%%%%%%%%%%%%%%%%%%%%%%%%%%%%%%%%%%%%%%%%%%%%%

\vspace{0.1cm}
Furthermore, in \cite{au-ma2}, Auvray, Ma, and Marinescu have derived optimal estimates of  $\mathcal{C}^{n}$-norms of the Bergman kernel associated to tensor-powers of $\mathcal{L}$, both along the diagonal, and away from the diagonal. The estimates derived in \cite{au-ma} and \cite{au-ma2} also remain stable in covers of Riemann surfaces. 
%%%%%%%%%%%%%%%%%%%%%%%%%%%%%%%%%%%%%%%%%%%%%%%%%%%%%%%%%%%%%%%%%%%%%%%%%%%%%%%%%%%%%%%%%%%%

\vspace{0.1cm}
We now state the main results from \cite{am}. Let $X$ be a hyperbolic Riemann surface of finite volume, and let $\O$ be the cotangent bundle of $X$, and let $\bkx$ denote the Bergman kernel associated to $\Ok$. Furthermore, let $\|\cdot\|_{\mathrm{hyp}}$ denote the point-wise hyperbolic metric on $\Ok$. Let $z,w\in X$ with $\delta:=\dh(z,w)\geq\rx$, where $\dh(z,w)$ denotes the geodesic distance between the points $z$ and $w$ on $X$, and $\rx$ is the injectivtiy radius of $X$, which is as defined in \eqref{irad1} or \eqref{irad2}, depending on  whether $X$ is compact or noncompact, respectively. 
%%%%%%%%%%%%%%%%%%%%%%%%%%%%%%%%%%%%%%%%%%%%%%%%%%%%%%%%%%%%%%%%%%%%%%%%%%%%%%%%%%%%%%%%%%%%

\vspace{0.1cm}
With hypothesis as above, when $X$ is  compact, we have the following estimate
\begin{align}\label{am1}
\|\bkx\|_{\mathrm{hyp}}(z,w)= O_{X}\left(\frac{k}{\cosh^{2k-4}\big((\delta-\rx)\slash 2\big)}\right).
\end{align} 
For $\delta \gg 0$, the above estimate is stronger than the estimate derived in \cite{ma} , which was derived in a more general context. We refer the reader to Remark $3.1$ in \cite{am}, for further  discussion on the above estimate.    
%%%%%%%%%%%%%%%%%%%%%%%%%%%%%%%%%%%%%%%%%%%%%%%%%%%%%%%%%%%%%%%%%%%%%%%%%%%%%%%%%%%%%%%%%%%%%

\vspace{0.1cm}
We then extend the above estimate to the noncompact setting, and show that the estimates remain stable in covers of Riemann surfaces. 
%%%%%%%%%%%%%%%%%%%%%%%%%%%%%%%%%%%%%%%%%%%%%%%%%%%%%%%%%%%%%%%%%%%%%%%%%%%%%%%%%%%%%%%%%%%%%

\vspace{0.1cm}\noindent
{\bf{Statement of Main theorem.}}
We now state the main theorem of the article.
%%%%%%%%%%%%%%%%%%%%%%%%%%%%%%%%%%%%%%%%%%%%%%%%%%%%%%%%%%%%%%%%%%%%%%%%%%%%%%%%%%%%%%%%%%%%%
\begin{theorem*}
With notation as above, for any $k\geq 3$, and $z=x+iy, \,w=u+iv\in X$ (identifying $X$ with its universal cover $\H$) with $\dh(z,w)=\delta$. Then, for $\rx\slash 2<\delta  <\rx$, when $X$ is compact, we have the following estimate
\begin{align}\label{mainthmeqn1}
&\|\bkx\|_{\mathrm{hyp}}(z,w)\leq \mathcal{C}_{X}^{1};
\end{align}
%%%%%%%%%%%%%%%%%%%%%%%%%%%%%%%%%%%%%%%%%%%%%%%%%%%%%%%%%%%%%%%%%%%%%%%%%%%%%%%%%%%%%%%%%%%
and when $X$ is noncompact, we have the following estimate
\begin{align}
&\|\bkx\|_{\mathrm{hyp}}(z,w)\leq \mathcal{C}_{X}^{1}+\frac{2k-1}{4\pi\cosh^{2k}\big(\delta\slash 2\big)}+\frac{(4yv)^{k}}{(y+v)^{2k-1}}\cdot\frac{(2k-1)\Gamma\big(k-1\slash 2\big)}{2\sqrt{\pi}\Gamma(k)},\label{mainthmeqn12}\\[0.1cm]
&\mathrm{where}\,\,\mathcal{C}_{X}^{1}:=\frac{2k-1}{4\pi}\bigg(\frac{1}{\cosh^{2k}\big((\rx-\delta)\slash 2\big)}+\frac{32}{\cosh^{2k-4}\big(\rx\slash 4\big)}\bigg)+\notag\\[0.1cm]&\hspace{3.7cm}\frac{2k-1}{\pi(k-2)\sinh^{2}\big(\rx\slash 4\big)}\cdot\frac{1}{\cosh^{2k-4}\big(\rx\slash 4\big)}.\label{constant1}
\end{align}
%%%%%%%%%%%%%%%%%%%%%%%%%%%%%%%%%%%%%%%%%%%%%%%%%%%%%%%%%%%%%%%%%%%%%%%%%%%%%%%%%%%%%%%%%%%
For $0\leq \delta \leq \rx\slash2$, when $X$ is compact, we have the following estimate
 \begin{align}\label{mainthmeqn2}
\|\bkx\|_{\mathrm{hyp}}(z,w)\leq \mathcal{C}_{X}^{2};
\end{align}
%%%%%%%%%%%%%%%%%%%%%%%%%%%%%%%%%%%%%%%%%%%%%%%%%%%%%%%%%%%%%%%%%%%%%%%%%%%%%%%%%%%%%%%%%%%
and when $X$ is noncompact, without loss of generality, we assume that $i\infty$ (identifying $X$ with its universal cover $\H$) is the only puncture of $X$. Then,  we have the following estimate
\begin{align}
&\|\bkx\|_{\mathrm{hyp}}(z,w)\leq \mathcal{C}_{X}^{2}+\frac{2k-1}{4\pi\cosh^{2k}\big(\delta\slash 2\big)}+\frac{(4yv)^{k}}{(y+v)^{2k-1}}\cdot\frac{(2k-1)\Gamma\big(k-1\slash 2\big)}{2\sqrt{\pi}\Gamma(k)},\label{mainthmeqn22}\\[0.1cm]
&\mathrm{where}\,\,\mathcal{C}_{X}^{2}:=
\frac{2k-1}{4\pi}\bigg(\frac{2}{\cosh^{2k}\big(\delta\slash 2\big)}+\frac{16}{\cosh^{2k-4}\big(\rx\slash 4\big)}+\frac{8}{\cosh^{2k-3}\big(\rx\slash 2\big)} \bigg)+
 \notag \\[0.1cm]&\frac{2k-1}{2\pi\sinh^{2}\big(\rx\slash4\big)}\cdot\bigg(\frac{1}{(2k-2)\cosh^{2k-3}\big(\rx\slash 2\big)}
+\frac{1}{(k-2)\cosh^{2k-4}\big(\rx\slash2\big)}\bigg).\label{constant2}
\end{align}
\end{theorem*}
%%%%%%%%%%%%%%%%%%%%%%%%%%%%%%%%%%%%%%%%%%%%%%%%%%%%%%%%%%%%%%%%%%%%%%%%%%%%%%%%%%%%%%%%%%%
%%%%%%%%%%%%%%%%%%%%%%%%%%%%%%%%%%%%%%%%%%%%%%%%%%%%%%%%%%%%%%%%%%%%%%%%%%%%%%%%%%%%%%%%%%%%%
%%%%%%%%%%%%%%%%%%%%%%%%%%%%%%%%%%%%%%%%%%%%%%%%%%%%%%%%%%%%%%%%%%%%%%%%%%%%%%%%%%%%%%%%%%%%%
%%%%%%%%%%%%%%%%%%%%%%%%%%%%%%%%%%%%%%%%%%%%%%%%%%%%%%%%%%%%%%%%%%%%%%%%%%%%%%%%%%%%%%%%%%%%%
\section{Background material}
%%%%%%%%%%%%%%%%%%%%%%%%%%%%%%%%%%%%%%%%%%%%%%%%%%%%%%%%%%%%%%%%%%%%%%%%%%%%%%%%%%%%%%%%%%%
We refer the reader to the section on background material in \cite{am}, for an elaborate discussion on the notation. However, we briefly explain the notation, and recall the results required for the proof of  Main theorem. 
%%%%%%%%%%%%%%%%%%%%%%%%%%%%%%%%%%%%%%%%%%%%%%%%%%%%%%%%%%%%%%%%%%%%%%%%%%%%%%%%%%%%%%%%%%%%%

\vspace{0.1cm}
Let $X$ be a hyperbolic Riemann surface of finite volume, which can be realized as the quotient space $\G\backslash\H$, where $\G\subset\mathrm{PSL}_{2}(\R)$ is a cofinite Fuchsian subgroup, and $\H$ is the complex upper half-plane. Locally, we identify $X$ with its universal cover $\H$, and hence, for only brevity of notation, we denote the points on $X$ by the same letters as the points on $\H$.
 %%%%%%%%%%%%%%%%%%%%%%%%%%%%%%%%%%%%%%%%%%%%%%%%%%%%%%%%%%%%%%%%%%%%%%%%%%%%%%%%%%%%%%%%%%%%%

\vspace{0.1cm}
Let $\dh(z,w)$ denote the hyperbolic distance on $\H$, which is the natural distance function on $\H$. Locally, for any $z,w\in X$, the geodesic distance between the points $z$ and $w$  on $X$ is given by $\dh(z,w)$. 
 %%%%%%%%%%%%%%%%%%%%%%%%%%%%%%%%%%%%%%%%%%%%%%%%%%%%%%%%%%%%%%%%%%%%%%%%%%%%%%%%%%%%%%%%%%%%

\vspace{0.1cm}
When $X$ is compact, injectivity radius $\rx$ is given by the following formula
\begin{align}\label{irad1}
\rx:=\inf\big\lbrace \dh(z,\gamma z)|\,z\in \H,\, \gamma\in\Gamma\backslash\lbrace\mathrm{Id}\rbrace\big\rbrace;
\end{align}
and when $X$ is noncompact, it is given by the following formula
\begin{align}\label{irad2}
\rx:=\inf\big\lbrace \dh(z,\gamma z)|\,z\in \H,\, \gamma\in\Gamma\backslash\G_{i\infty}\big\rbrace,
\end{align}
where $\G_{i\infty}$ is the stabilizer of $i\infty$. Here, as stated in Main theorem, identifying $X$ with its universal cover $\H$, we assume that $i\infty$ is the only puncture of $X$.   
%%%%%%%%%%%%%%%%%%%%%%%%%%%%%%%%%%%%%%%%%%%%%%%%%%%%%%%%%%%%%%%%%%%%%%%%%%%%%%%%%%%%%%%%%%%%%

\vspace{0.1cm}
Let $\O$ denote the cotangent bundle of holomorphic differential $1$-forms on $X$. Then, for any $k>0$, the Bergman kernel associated to $\Ok$ is given by the following formula
\begin{align*}
&\bkx(z,w):= \frac{(2k-1)(2i)^{2k}}{4\pi}\sum_{\gamma\in\Gamma}\frac{1}{\big(\gamma z-\overline{w}\big)^{2k}}\cdot\frac{\big(dz^{\otimes k}\wedge d\overline{w}^{\otimes k}\big)}{j\big(\gamma,z\big)^{2k}}\\
&\mathrm{where},\,\,\mathrm{for\,\,any}\,\gamma=\left(\begin{array}{cc} a&b\\c&d\end{array}\right)\in\Gamma, \,j\big(\gamma,z\big)=cz+d.\,
\end{align*}
%%%%%%%%%%%%%%%%%%%%%%%%%%%%%%%%%%%%%%%%%%%%%%%%%%%%%%%%%%%%%%%%%%%%%%%%%%%%%%%%%%%%%%%%%%%%%

\vspace{0.025cm}
The hyperbolic metric on $\Ok$ induces the following point-wise hyperbolic metric on $\bkx(z,w)$
\begin{align*}
\|\bkx\|_{\mathrm{hyp}}(z,w)= \frac{(2k-1)(4yv)^{k}}{4\pi}\cdot\Bigg|\sum_{\gamma\in\Gamma}\frac{1}{\big(\gamma z-\overline{w}\big)^{2k}}\cdot\frac{1}{j(\gamma,z)^{2k}}\Bigg|.
\end{align*}
%%%%%%%%%%%%%%%%%%%%%%%%%%%%%%%%%%%%%%%%%%%%%%%%%%%%%%%%%%%%%%%%%%%%%%%%%%%%%%%%%%%%%%%%%%%%%
For any $z = x + iy, \,w = u + iv \in \H$, and for any $\gamma\in \Gamma$, we have the following two formulae
\begin{align}
\mathrm{Im}\big(\gamma z\big)=\frac{y}{\big|cz+d\big|^{2}}\,\,\mathrm{and}\,\,\cosh^{2}\big(\dh(z,w)\slash2\big)= \frac{|z -w|^{2}}{4yv}.
 \end{align}
 %%%%%%%%%%%%%%%%%%%%%%%%%%%%%%%%%%%%%%%%%%%%%%%%%%%%%%%%%%%%%%%%%%%%%%%%%%%%%%%%%%%%%%%%%%%%%

\vspace{0.025cm}
For any $z,w\in \H$, combining the above two equations, we derive the following inequality
\begin{align}\label{seriesbergmanreln}
\|\bkx\|_{\mathrm{hyp}}(z,w)\leq \frac{2k-1}{4\pi}\sum_{\gamma\in\Gamma}\frac{\big(4\mathrm{Im}\big(\gamma z\big)\cdot v\big)^{k}}{\big|\gamma z-\overline{w}\big|^{2k}}=
\frac{2k-1}{4\pi}\sum_{\gamma\in\Gamma}\frac{1}{\cosh^{2k}\big(\dh(\gamma z,w)\slash 2\big)}.
\end{align}
%%%%%%%%%%%%%%%%%%%%%%%%%%%%%%%%%%%%%%%%%%%%%%%%%%%%%%%%%%%%%%%%%%%%%%%%%%%%%%%%%%%%%%%%%%%%%

\vspace{0.1cm}
We now state two inequalities from \cite{jl}, which are adapted to our setting. The inequalities give us an estimate for the number of elements in $\G$ or $\G\backslash \G_{i\infty}$, depending on whether $X$ is compact or noncompact, respectively.
%%%%%%%%%%%%%%%%%%%%%%%%%%%%%%%%%%%%%%%%%%%%%%%%%%%%%%%%%%%%%%%%%%%%%%%%%%%%%%%%%%%%%%%%%%%%%

For any positive, smooth, real-valued, and decreasing function $f$ defined on $\mathbb{R}_{\geq 0}$, and for any $\delta > \rx\slash 2$, and $z,w\in \H$, we have the following inequality 
\begin{align}\label{jlineq1}
 &\int_0^{\infty} f(\rho) dN_{\G}(z,w;\rho) \leq \int_0^{\delta} f(\rho) dN_{\Gamma}\big(z,w;\rho\big) + f(\delta) \frac{2\cosh\big(\rx\slash 4\big)\sinh(\delta)}{\sinh\big(\rx\slash 4\big)}
 +\notag\\[0.1cm]&\hspace{5.4cm} \frac{1}{2\sinh^2\big(\rx\slash 4\big)} \int_{\delta}^{\infty} f(\rho) \sinh\big(\rho +\rx\slash 2\big) d\rho;
\end{align}  
%%%%%%%%%%%%%%%%%%%%%%%%%%%%%%%%%%%%%%%%%%%%%%%%%%%%%%%%%%%%%%%%%%%%%%%%%%%%%%%%%%%%%%%%%%%%%
for any $\delta>0$, we have the following inequality
\begin{align}\label{jlineq2}
N_{\G}\big(z,w;\delta\big)\leq \frac{\sinh\big(\delta+\rx\big)}{\sinh\big(\rx\big)}, \,\, \mathrm{where}\,\,N_{\Gamma}\big(z,w;\rho\big) := \mathrm{card}\,\{ \gamma |\,\gamma \in \G\backslash\Gamma_{\infty},\,\dh(\gamma z,w) \leq \rho \}.
\end{align}
%%%%%%%%%%%%%%%%%%%%%%%%%%%%%%%%%%%%%%%%%%%%%%%%%%%%%%%%%%%%%%%%%%%%%%%%%%%%%%%%%%%%%%%%%%%%%

In the above inequalities, it is implicitly understood that, when $X$ is compact, $\G_{i\infty}=\emptyset$. 
%%%%%%%%%%%%%%%%%%%%%%%%%%%%%%%%%%%%%%%%%%%%%%%%%%%%%%%%%%%%%%%%%%%%%%%%%%%%%%%%%%%%%%%%%%%%%
%%%%%%%%%%%%%%%%%%%%%%%%%%%%%%%%%%%%%%%%%%%%%%%%%%%%%%%%%%%%%%%%%%%%%%%%%%%%%%%%%%%%%%%%%%%%%
%%%%%%%%%%%%%%%%%%%%%%%%%%%%%%%%%%%%%%%%%%%%%%%%%%%%%%%%%%%%%%%%%%%%%%%%%%%%%%%%%%%%%%%%%%%%%
\section{Proof of the Main theorem}
%%%%%%%%%%%%%%%%%%%%%%%%%%%%%%%%%%%%%%%%%%%%%%%%%%%%%%%%%%%%%%%%%%%%%%%%%%%%%%%%%%%%%%%%%%%%%
\begin{proof}[Proof of estimate \eqref{mainthmeqn1}]
For any $k\geq 3$, and $z,w\in X$ with $\rx\slash 2<\dh(z,w) = \delta<\rx$, using inequalities \eqref{seriesbergmanreln} and \eqref{jlineq1}, we find that
\begin{align}
&\|\bkx\|_{\mathrm{hyp}}(z,w) \leq \frac{2k-1}{4\pi}\sum_{\gamma\in\Gamma}\frac{1}{\cosh^{2k}\big(\dh(\gamma z,w)\slash 2\big)}=\notag\\[0.1cm]
&\frac{2k-1}{4\pi}\int_{0}^{\infty}\frac{dN_{\G}\big(z,w;\rho\big)}{\cosh^{2k}\big(\dh(\gamma z,w)\slash 2\big)}\leq \frac{2k-1}{4\pi}\int_{0}^{\delta}\frac{dN_{\G}\big(z,w;\rho\big)}{\cosh^{2k}\big(\dh(\gamma z,w)\slash 2\big)}+\notag\\[0.1cm]
&\frac{2k-1}{4\pi}\cdot\frac{2\cosh\big(\rx\slash4\big)\sinh(\delta)}{\cosh^{2k}\big(\delta\slash2\big)\sinh\big(\rx\slash 4\big)}+\frac{2k-1}{8\pi\sinh^{2}\big(\rx\slash 4\big)}\int_{\delta}^{\infty}\frac{\sinh\big(\rho+\rx\slash 2\big)d\rho}{\cosh^{2k}\big(\rho\slash 2\big)}.\label{proof1eqn1}
\end{align}
%%%%%%%%%%%%%%%%%%%%%%%%%%%%%%%%%%%%%%%%%%%%%%%%%%%%%%%%%%%%%%%%%%%%%%%%%%%%%%%%%%%%%%%%%%%%%

\vspace{0.1cm}
We now estimate the first term on the right hand-side of the above inequality
\begin{align*}
&\int_{0}^{\delta}\frac{dN_{\G}(z,w;\rho)}{\cosh^{2k}\big(\dh(\gamma z,w)\slash 2\big)}=
\int_{0}^{\rx\slash2}\frac{dN_{\G}(z,w;\rho)}{\cosh^{2k}\big(\dh(\gamma z,w)\slash 2\big)}+\int_{\rx\slash 2}^{\delta}\frac{dN_{\G}(z,w;\rho)}{\cosh^{2k}\big(\dh(\gamma z,w)\slash 2\big)}.
\end{align*}
%%%%%%%%%%%%%%%%%%%%%%%%%%%%%%%%%%%%%%%%%%%%%%%%%%%%%%%%%%%%%%%%%%%%%%%%%%%%%%%%%%%%%%%%%%%%%

\vspace{0.1cm}
From the definition of injectivity radius (equation \eqref{irad1}), it is clear that there can be at most one $\gamma\in \Gamma$ such that $\dh(\gamma z,w)<\rx\slash 2$. Furthermore, for any $\gamma\in \Gamma$, using triangular inequality, we observe that 
\begin{align*}
\dh(\gamma z,w)+\dh(w,z)\geq \dh(z,\gamma z)&\geq \rx\implies \dh(\gamma z,w)\geq\rx-\delta\implies \notag\\\frac{1}{\cosh^{2k}\big(\dh(\gamma z,w)\slash 2\big)}&\leq \frac{1}{\cosh^{2k}\big((\rx-\delta)\slash 2\big)}.
\end{align*}
  %%%%%%%%%%%%%%%%%%%%%%%%%%%%%%%%%%%%%%%%%%%%%%%%%%%%%%%%%%%%%%%%%%%%%%%%%%%%%%%%%%%%%%%%%%%%%

\vspace{0.1cm}
From the above observations, we arrive at the following inequality
 \begin{align}\label{proof1eqn2}
 \int_{0}^{\rx\slash 2}\frac{dN_{\G}(z,w;\rho)}{\cosh^{2k}\big(\dh(\gamma z,w)\slash 2\big)}\leq \frac{1}{\cosh^{2k}\big((\rx-\delta)\slash 2\big)}.
 \end{align}
  %%%%%%%%%%%%%%%%%%%%%%%%%%%%%%%%%%%%%%%%%%%%%%%%%%%%%%%%%%%%%%%%%%%%%%%%%%%%%%%%%%%%%%%%%%%%%  

\vspace{0.1cm}
Using inequality \eqref{jlineq2}, and the hypothesis that $\delta<\rx$, and combining it with the observation that $\cosh(x)\leq 2\cosh^{2}(x\slash 2)$, for all $x\geq 0$, we derive the following inequality
\begin{align}\label{proof1eqn3}
&\int_{\rx\slash 2}^{\delta}\frac{dN_{\G}\big(z,w;\rho\big)}{\cosh^{2k}\big(\dh(\gamma z,w)\slash 2\big)}\leq \sup_{\rho\in[\rx\slash 2,\delta]}\frac{
N_{\G}\big(z,w;\delta\big)}{\cosh^{2k}\big(\rho\slash 2\big)}\leq\notag\\[0.1cm] &\frac{\sinh\big(2\rx\big)}{\cosh^{2k}\big(\rx\slash4\big)\sinh\big(\rx\big)}=\frac{2\cosh\big(\rx\big)}{\cosh^{2k}\big(\rx\slash4\big)}\leq\frac{16}{\cosh^{2k-4}\big(\rx\slash4\big)}.
\end{align}
%%%%%%%%%%%%%%%%%%%%%%%%%%%%%%%%%%%%%%%%%%%%%%%%%%%%%%%%%%%%%%%%%%%%%%%%%%%%%%%%%%%%%%%%%%%%%

\vspace{0.1cm}
Using the hypothesis that $\rx\slash 2<\delta<\rx$, and combining it with the observation that $\cosh(x)\leq 2\cosh^{2}(x\slash 2)$, for all $x\geq 0$, we arrive at the following estimate for the second term on the right hand-side of inequality \eqref{proof1eqn1}
\begin{align}\label{proof1eqn4}
\frac{2\cosh\big(\rx\slash4\big)\sinh(\delta)}{\cosh^{2k}\big(\delta\slash2\big)\sinh\big(\rx\slash 4\big)}\leq \frac{8\cosh\big(\rx\slash2\big)}{\cosh^{2k-2}\big(\rx\slash4\big)}\leq \frac{16}{\cosh^{2k-4}\big(\rx\slash4\big)}.
\end{align}
%%%%%%%%%%%%%%%%%%%%%%%%%%%%%%%%%%%%%%%%%%%%%%%%%%%%%%%%%%%%%%%%%%%%%%%%%%%%%%%%%%%%%%%%%%%%%

\vspace{0.1cm}
We have the following inequality from inequality $(12)$ in \cite{am}
\begin{align}\label{proof1eqn5}
\int_{\delta}^{\infty}\frac{\sinh\big(\rho+\rx\slash 2\big)d\rho}{\cosh^{2k}\big(\rho\slash 2\big)}\leq \frac{4\cosh\big(\rx\slash 2\big)}{(2k-2) \cosh^{2k-2}\big(\delta\slash2\big)} + \frac{8}{(2k-4)\cosh^{2k-4}\big(\delta\slash2\big)}.
\end{align}
%%%%%%%%%%%%%%%%%%%%%%%%%%%%%%%%%%%%%%%%%%%%%%%%%%%%%%%%%%%%%%%%%%%%%%%%%%%%%%%%%%%%%%%%%%%%%

\vspace{0.1cm}
Using the above inequality, and the hypothesis that $\rx\slash 2<\delta<\rx$, and combining it with the observation that $\cosh(x)\leq 2\cosh^{2}(x\slash 2)$, for all $x\geq 0$, we derive that
\begin{align}\label{proof1eqn6}
&\int_{\delta}^{\infty}\frac{\sinh\big(\rho+\rx\slash 2\big)d\rho}{\cosh^{2k}\big(\rho\slash 2\big)}\leq\notag\\[0.1cm]&  \frac{8}{(2k-2)\cosh^{2k-4}\big(\rx\slash4\big)}+\frac{8}{(2k-4)\cosh^{2k-4}\big(\rx\slash4\big)}\leq \frac{16}{(2k-4)\cosh^{2k-4}\big(\rx\slash4\big)}.
\end{align}
%%%%%%%%%%%%%%%%%%%%%%%%%%%%%%%%%%%%%%%%%%%%%%%%%%%%%%%%%%%%%%%%%%%%%%%%%%%%%%%%%%%%%%%%%%%%%

\vspace{0.1cm}
Combining estimates \eqref{proof1eqn1}, \eqref{proof1eqn2}, \eqref{proof1eqn3}, \eqref{proof1eqn4}, and \eqref{proof1eqn6} completes the proof of estimate \eqref{mainthmeqn1}.
\end{proof}
%%%%%%%%%%%%%%%%%%%%%%%%%%%%%%%%%%%%%%%%%%%%%%%%%%%%%%%%%%%%%%%%%%%%%%%%%%%%%%%%%%%%%%%%%%%%%
%%%%%%%%%%%%%%%%%%%%%%%%%%%%%%%%%%%%%%%%%%%%%%%%%%%%%%%%%%%%%%%%%%%%%%%%%%%%%%%%%%%%%%%%%%%%%
%%%%%%%%%%%%%%%%%%%%%%%%%%%%%%%%%%%%%%%%%%%%%%%%%%%%%%%%%%%%%%%%%%%%%%%%%%%%%%%%%%%%%%%%%%%%%

\vspace{0.2cm}
\begin{proof}[Proof of estimate \eqref{mainthmeqn2}]
For any $k\geq 3$, and $z,w\in X$ with $0\leq \dh(z,w) = \delta\leq \rx\slash 2$, using \eqref{jlineq1}, we find that
\begin{align}\label{proof2eqn1}
&\|\bkx\|_{\mathrm{hyp}}(z,w) \leq \frac{2k-1}{4\pi}\int_{0}^{\rx}\frac{dN_{\G}\big(z,w;\rho\big)}{\cosh^{2k}\big(\dh(\gamma z,w)\slash 2\big)}
+\notag\\[0.1cm]&\frac{2k-1}{4\pi\cosh^{2k}\big(\rx\slash 2\big)}\cdot\frac{2\cosh\big(\rx\slash 4\big)\sinh\big(\rx\big)}{\sinh\big(\rx\slash 4\big)}+\frac{2k-1}{8\pi\sinh^{2}\big(\rx\slash4\big)}\int_{\rx}^{\infty}\frac{\sinh\big(\rho +\rx\slash 2\big) d\rho}{\cosh^{2k}\big(\rho\slash 2\big)}.
\end{align}
%%%%%%%%%%%%%%%%%%%%%%%%%%%%%%%%%%%%%%%%%%%%%%%%%%%%%%%%%%%%%%%%%%%%%%%%%%%%%%%%%%%%%%%%%%%%%

\vspace{0.1cm}
We now estimate the first term on the right hand side of the above inequality
\begin{align*}
\int_{0}^{\rx}\frac{dN_{\G}\big(z,w;\rho\big)}{\cosh^{2k}\big(\dh(\gamma z,w)\slash 2\big)}=\int_{0}^{\rx\slash 2}\frac{dN_{\G}\big(z,w;\rho\big)}{\cosh^{2k}\big(\dh(\gamma z,w)\slash 2\big)}+\int_{\rx\slash 2}^{\rx}\frac{dN_{\G}\big(z,w;\rho\big)}{\cosh^{2k}\big(\dh(\gamma z,w)\slash 2\big).}
\end{align*}
%%%%%%%%%%%%%%%%%%%%%%%%%%%%%%%%%%%%%%%%%%%%%%%%%%%%%%%%%%%%%%%%%%%%%%%%%%%%%%%%%%%%%%%%%%%%%

\vspace{0.1cm}
For $0\leq \dh(z,w) = \delta\leq \rx\slash 2$, from the definition of injectivity radius (equation \eqref{irad1}), there can be at most one more $\gamma\in \Gamma$ other than 
$\gamma={\rm{Id}}$ such that $\dh(\gamma z,w)\leq \rx\slash 2$. So we have
\begin{align}\label{proof2eqn2}
\int_{0}^{\rx\slash 2}\frac{dN_{\G}\big(z,w;\rho\big)}{\cosh^{2k}\big(\dh(\gamma z,w)\slash 2\big)}\leq \frac{2}{\cosh^{2k}\big(\delta\slash 2\big)}.
\end{align}
%%%%%%%%%%%%%%%%%%%%%%%%%%%%%%%%%%%%%%%%%%%%%%%%%%%%%%%%%%%%%%%%%%%%%%%%%%%%%%%%%%%%%%%%%%%%%

\vspace{0.025cm}
From similar arguments as in inequality \eqref{proof1eqn3}, we arrive at the following inequality
\begin{align}\label{proof2eqn3}
&\int_{\rx\slash 2}^{\rx}\frac{dN_{\G}\big(z,w;\rho\big)}{\cosh^{2k}\big(\dh(\gamma z,w)\slash 2\big)}\leq \frac{N_{\G}\big(z,w;\rx\big)}{\cosh^{2k}\big(\rx\slash 4\big)}\leq  \frac{2\cosh\big(\rx\big)}{\cosh^{2k}\big(\rx\slash 4\big)}\leq\frac{16}{\cosh^{2k-4}\big(\rx\slash 4\big)}.
\end{align}
%%%%%%%%%%%%%%%%%%%%%%%%%%%%%%%%%%%%%%%%%%%%%%%%%%%%%%%%%%%%%%%%%%%%%%%%%%%%%%%%%%%%%%%%%%%%%

\vspace{0.1cm}
Using similar arguments as in the proof of estimate \eqref{mainthmeqn2}, we derive the following inequality
\begin{align}\label{proof2eqn4}
&\frac{2}{\cosh^{2k}\big(\rx\slash 2\big)}\cdot\frac{\cosh\big(\rx\slash 4\big)\sinh\big(\rx\big)}{\sinh\big(\rx\slash 4\big)}\leq \frac{8}{\cosh^{2k-3}\big(\rx\slash 2\big)}.
\end{align}
%%%%%%%%%%%%%%%%%%%%%%%%%%%%%%%%%%%%%%%%%%%%%%%%%%%%%%%%%%%%%%%%%%%%%%%%%%%%%%%%%%%%%%%%%%%%%

\vspace{0.1cm}
Substituting $\delta=\rx$ in inequality \eqref{proof1eqn5}, we arrive at the following inequality  
\begin{align}
&\int_{\rx}^{\infty}\frac{\sinh\big(\rho+\rx\slash 2\big)d\rho}{\cosh^{2k}\big(\rho\slash 2\big)}\leq \frac{4}{(2k-2)\cosh^{2k-3}\big(\rx\slash 2\big)}+\frac{8}{(2k-4)\cosh^{2k-4}\big(\rx\slash 2\big)}. \label{proof2eqn5}
\end{align}
%%%%%%%%%%%%%%%%%%%%%%%%%%%%%%%%%%%%%%%%%%%%%%%%%%%%%%%%%%%%%%%%%%%%%%%%%%%%%%%%%%%%%%%%%%%%%

\vspace{0.1cm}
Combining estimates \eqref{proof2eqn1}, \eqref{proof2eqn2}, \eqref{proof2eqn3}, \eqref{proof2eqn4}, and \eqref{proof2eqn5} completes the proof of estimate \eqref{mainthmeqn2}.
\end{proof}
%%%%%%%%%%%%%%%%%%%%%%%%%%%%%%%%%%%%%%%%%%%%%%%%%%%%%%%%%%%%%%%%%%%%%%%%%%%%%%%%%%%%%%%%%%%%%
%%%%%%%%%%%%%%%%%%%%%%%%%%%%%%%%%%%%%%%%%%%%%%%%%%%%%%%%%%%%%%%%%%%%%%%%%%%%%%%%%%%%%%%%%%%%%

\vspace{0.2cm}
\begin{proof}[Proofs of estimates \eqref{mainthmeqn12} and \eqref{mainthmeqn22}]
The proofs of estimates \eqref{mainthmeqn12} and \eqref{mainthmeqn22} follow from the proof of estimate $(2)$ from \cite{am}, when combined with the proofs of estimates \eqref{mainthmeqn1} and \eqref{mainthmeqn2}, respectively. 
\end{proof}
%%%%%%%%%%%%%%%%%%%%%%%%%%%%%%%%%%%%%%%%%%%%%%%%%%%%%%%%%%%%%%%%%%%%%%%%%%%
%%%%%%%%%%%%%%%%%%%%%%%%%%%%%%%%%%%%%%%%%%%%%%%%%%%%%%%%%%%%%%%%%%%%%%%%%%%%%%%%%%%%%%%%%%%%%

\vspace{0.1cm}
\begin{rem}
With hypothesis as in Main theorem, let $X$ be a compact hyperbolic Riemann surface, and let $z,w\in X$ with $\delta:=\dh(z,w)$. If $\rx\slash 2<\delta<\rx$, a careful analysis of each of the term comprising estimate \eqref{mainthmeqn1} leads us to the following estimate
\begin{align}
&\|\bkx\|_{\mathrm{hyp}}(z,w)=O_{X}\left(\frac{k}{\cosh^{2k}\big((\rx-\delta)\slash 2\big)}\right);\label{remeqn1}
&\intertext{for any $0\leq \delta\leq \rx\slash 2$, we have the following estimate}
&\|\bkx\|_{\mathrm{hyp}}(z,w)=O_{X}\left(\frac{k}{\cosh^{2k}\big(\delta\slash 2\big)}\right).\label{remeqn2}
\end{align}
\end{rem}
%%%%%%%%%%%%%%%%%%%%%%%%%%%%%%%%%%%%%%%%%%%%%%%%%%%%%%%%%%%%%%%%%%%%%%%%%%%%%%%%%%%%%%%%%%%%%

\vspace{0.1cm}
\begin{rem}
With hypothesis as in Main theorem, let $X$ be a compact hyperbolic Riemann surface. Along the diagonal, when $z=w\in X$, from Main theorem, it is easy to derive the following estimate 
\begin{align}\label{remeqn3}
\|\bkx\|_{\mathrm{hyp}}(z,z)=O_{X}(k).
\end{align}
%%%%%%%%%%%%%%%%%%%%%%%%%%%%%%%%%%%%%%%%%%%%%%%%%%%%%%%%%%%%%%%%%%%%%%%%%%%%%%%%%%%%%%%%%%%%%

\vspace{0.1cm}
When $X$ is noncompact, from the proof of Proposition $5.1$ in $p. 11$ and $p. 12$ in \cite{jk}, it is clear that the Bergman kernel
$\|\bkx\|_{\mathrm{hyp}}(z,z)$ takes its maximum value on $\partial D$, which is the boundary of the following strip 
\begin{align*}
\mathcal{D}:=\bigg\lbrace z=x+iy\in\H |\,0\leq x\leq 1,\,y> \frac{k}{2\pi} \bigg\rbrace ,
\end{align*}
%%%%%%%%%%%%%%%%%%%%%%%%%%%%%%%%%%%%%%%%%%%%%%%%%%%%%%%%%%%%%%%%%%%%%%%%%%%%%%%%%%%%%%%%%%%%%
which implies that 
\begin{align}\label{rem2eqn1}
&\sup_{z\in\H}\|\bkx\|_{\mathrm{hyp}}(z,z)\leq\notag\\&\sup_{z\in\partial \mathcal{D}}\frac{2k-1}{4\pi} \sum_{\gamma\in\Gamma\backslash\Gamma_{\infty}} \frac{1}{\cosh^{2k}\big(\dh(\gamma z,z)\slash 2\big)}+\sup_{z\in\partial \mathcal{D}}\frac{2k-1}{4\pi} \sum_{\gamma\in\Gamma_{\infty}} \frac{1}{\cosh^{2k}\big(\dh(\gamma z,z)\slash 2\big)}.
\end{align}  
%%%%%%%%%%%%%%%%%%%%%%%%%%%%%%%%%%%%%%%%%%%%%%%%%%%%%%%%%%%%%%%%%%%%%%%%%%%%%%%%%%%%%%%%%%%%%

\vspace{0.1cm}
From the arguments from the proof of Main theorem, it is clear that the first term on the right hand side of the above inequality satisfies the following estimate
\begin{align}\label{rem2eqn2}
\sup_{z\in\partial \mathcal{D}}\frac{2k-1}{4\pi} \sum_{\gamma\in\Gamma\backslash\Gamma_{\infty}} \frac{1}{\cosh^{2k}\big(\dh(\gamma z,z)\slash 2\big)}=O_{X}(k).
\end{align}
%%%%%%%%%%%%%%%%%%%%%%%%%%%%%%%%%%%%%%%%%%%%%%%%%%%%%%%%%%%%%%%%%%%%%%%%%%%%%%%%%%%%%%%%%%%%%

\vspace{0.1cm}
From the arguments from the proof of estimate $(2)$ in \cite{am}, and from asymptotics of the Gamma function, we have the following estimate for the second term on the right hand-side of inequality \eqref{rem2eqn1}
\begin{align}\label{rem2eqn3}
&\sup_{z\in\partial \mathcal{D}}\frac{2k-1}{4\pi} \sum_{\gamma\in\Gamma_{\infty}} \frac{1}{\cosh^{2k}\big(\dh(\gamma z,z)\slash 2\big)}\leq \notag\\&\frac{2k-1}{4\pi}+\sup_{z\in\partial \mathcal{D}} \frac{y\cdot(2k-1)\Gamma\big(k-1\slash 2\big)}{\sqrt{\pi}\Gamma(k)}= O\big(k^{3\slash 2}\big).
\end{align}
%%%%%%%%%%%%%%%%%%%%%%%%%%%%%%%%%%%%%%%%%%%%%%%%%%%%%%%%%%%%%%%%%%%%%%%%%%%%%%%%%%%%%%%%%%%%%

\vspace{0.1cm}
Therefore, for any $z=w\in X$, combining estimates \eqref{rem2eqn1}, \eqref{rem2eqn2}, and \eqref{rem2eqn3}, we arrive at the following estimate
\begin{align}\label{rem2eqn4}
\|\bkx\|_{\mathrm{hyp}}(z,z)=O\big(k^{3\slash 2}\big).
\end{align}
%%%%%%%%%%%%%%%%%%%%%%%%%%%%%%%%%%%%%%%%%%%%%%%%%%%%%%%%%%%%%%%%%%%%%%%%%%%%%%%%%%%%%%%%%%%%%
Estimates \eqref{remeqn3} and \eqref{rem2eqn4} have already been proved in \cite{jk} and \cite{au-ma2}, using a different approach, and were shown to be optimal. 
\end{rem}
%%%%%%%%%%%%%%%%%%%%%%%%%%%%%%%%%%%%%%%%%%%%%%%%%%%%%%%%%%%%%%%%%%%%%%%%%%%%%%%%%%%%%%%%%%%%%

\vspace{0.1cm}
\begin{rem}
From arguments similar to the ones employed in Remark 3.3 in \cite{am}, it is easy to show that estimates \eqref{remeqn1}, \eqref{remeqn2}, \eqref{remeqn3}, and \eqref{rem2eqn4} remain stable in covers of Riemann surfaces. 
\end{rem}
%%%%%%%%%%%%%%%%%%%%%%%%%%%%%%%%%%%%%%%%%%%%%%%%%%%%%%%%%%%%%%%%%%%%%%%%%%%%%%%%%%%%%%%%%%%%%
{\bf{Acknowledgements}}.  Both the authors acknowledge the support of INSPIRE research grant DST/INSPIRE/04/2015/002263.
%%%%%%%%%%%%%%%%%%%%%%%%%%%%%%%%%%%%%%%%%%%%%%%%%%%%%%%%%%%%%%%%%%%%%%%%%%%%%%%%%%%%%%%%%%%%%
%%%%%%%%%%%%%%%%%%%%%%%%%%%%%%%%%%%%%%%%%%%%%%%%%%%%%%%%%%%%%%%%%%%%%%%%%%%%%%%%%%%%%%%%%%%%%

\end{document}